# Relações Matemáticas: uma ferramenta no combate ao desinteresse dos alunos

# Mathematical Relationships: a toll in the fight against students' disinterest

Maurício Paulino Marques Fernandes[1]

**RESUMO**

Neste artigo, sugerimos o uso de Relações Matemáticas como um possível caminho para diminuir o desinteresse dos alunos pelas exatas. Primeiramente fizemos uma consideração a respeito do ambiente da sala de aula na perspectiva de professores e alunos e comparamos com o que é proposto nos PCN´s, procuramos mostrar que a aplicação dos conceitos matemáticos já faz parte da dia a dia dos alunos, embora eles não os relacionem aos conteúdos ministrados nas escolas, também trouxemos alguns exemplos de exploração de ambientes, situações e filmes que auxiliam na elaboração de temas geradores para as aulas de matemática, por fim, convidamos os professores a refletirem sobre seu papel como mediadores em sala de aula.

Palavras Chaves: Matemática Aplicada. Educação Matemática. Situações Problemas.

**Abstract**

In this article, we suggest the use of Mathematical Relationships as a possible way to decrease the students' disinterest in Mathematics. First we made a consideration of the environment of the classroom from the perspective of teachers and students and compare it with what is proposed in PCN's, sought to show that the application of mathematical concepts is already part of the daily lives of students, although they not relate to the material taught in schools, also brought some examples of operating environments, situations and movies that help in the preparation of generating themes for math classes, finally, we invite teachers to reflect on their roles as mediators in classroom.

---

[1] Pós Graduando em Matemática pela Universidade Estadual de Campinas – Unicamp, Campinas – SP, Brasil, E-mail: m130099@dac.unicamp.br

**Introdução**

Os recentes resultados, amplamente divulgados pela mídia em relação ao aprendizado matemático revelam a preocupante situação dos estudantes brasileiros; apesar de vermos grandes conquistas em competições nacionais e internacionais nas áreas das exatas, a realidade nos mostra que a maior parcela da sociedade não domina a matemática a contento do previsto nos PCN's, principalmente no que diz respeito ao conhecimento tecnológico esperado ao final do ensino médio.

> *"No nível médio, esses objetivos envolvem, de um lado, o aprofundamento dos saberes disciplinares em Biologia, Física, Química e Matemática, com procedimentos científicos pertinentes aos seus objetos de estudo, com metas formativas particulares, até mesmo com tratamentos didáticos específicos. De outro lado, envolvem a articulação interdisciplinar desses saberes, propiciada por várias circunstâncias, dentre as quais se destacam os conteúdos tecnológicos e práticos, já presentes junto a cada disciplina, mas particularmente apropriados para serem tratados desde uma perspectiva integradora."* (Parâmetros Curriculares Nacionais)

Dentre as várias causas apontadas para tal situação, uma das mais difundidas entre a sociedade e principalmente entre os professores é a que diz respeito ao desinteresse geral dos alunos, que por sua vez não concebem a escola como o espaço de formação já que ela, de fato, não cumpre seu papel de ligação entre currículo e realidade. Em matemática não é diferente, a dificuldade de relacioná-la com objetos e situações do cotidiano colabora para o desinteresse, dando a impressão de que o que se ensina não possui aplicação prática.

Para Skovsmose o aprendizado matemático deve contemplar três tipos de conhecimentos: Matemático, Tecnológico e Reflexivo. O conhecimento matemático envolve a teoria em si, as técnicas de cálculos e desenvolvimentos dedutivos em teoremas e suas demonstrações, o conhecimento tecnológico refere-se à aplicação do conhecimento matemático, é a construção de modelos matemáticos e deve ser aprendido na escola por meio de atividade de aplicação matemática em situações problemas. O conhecimento reflexivo é a habilidade de avaliar o uso da matemática, de refletir sobre tal uso, de fazer inferências nos métodos e resultados.

**Relações Matemáticas como Alternativa**

Negligenciar os conhecimentos tecnológico e reflexivo, e contemplar apenas o conhecimento matemático teórico, tem sido, na maioria das vezes, prática comum nas escolas brasileiras. Em contra partida vivemos em um mundo globalizado, de informações rápidas e precisas, aprender matemática ou qualquer outra ciência sem que ela seja precisamente útil, não faz sentido para os alunos dessa geração, é preciso levar os alunos à contextualizarem conceitos, pensar soluções e refletir sobre elas, esse é de fato o papel do professor nos dias atuais, não cabe mais a figura do professor que transmite informações e conteúdos, pois eles estão ao alcance das mãos e dos olhos, basta uma simples pesquisa na internet e pronto, toda a informação que buscamos surge como num passe de mágicas.

É preciso, portanto repensar o papel do professor, e esse repensar, não se restringe aos órgãos responsáveis pela educação ou às universidades, responsável pela formação de professores, cabe antes e principalmente ao professor, pois é ele quem estará a frente dos estudantes. Faz-se necessário considerar a tecnologia, os meios de comunicação e as redes sociais que nos cercam, aprender a lidar com essas ferramentas e utilizá-las de forma a construir caminhos que nos levem a uma relação mais próxima de cada aluno de forma que possam compreender o papel da escola e seus ensinos.

Quando falamos em aprendizado tecnológico, não nos referimos a um caminho árduo a percorrer, preferimos a palavra trabalhosa, pois o caminho não é sofrível, antes exige dedicação, treino e principalmente observação. O professor que pretende realmente diversificar suas aulas, contextualizá-la e ainda diminuir o desinteresse necessita estar de olhos bem abertos, pois oportunidades de contextualização aparecem em todos os lugares, quanto mais atento estiver, mais relações conseguirá formular para suas aulas. Sugerimos portanto a Matemática Aplicada como uma ferramenta no combate ao desinteresse dos alunos e como estímulo ao raciocínio lógico matemático.

Entendemos por Relações Matemáticas todo e qualquer conjunto de ações ou situações problemas que levem à flexão matemática, sejam estas ações e/ou reflexões simples ou não. Toda a reflexão é por si só um avanço no processo de aprendizagem e a medida em que se torna prática constante nas salas de aula, tornar-se-á uma atividade comum e simples tanto para alunos quanto professores, criando hábitos de observações no mundo atual, propiciando a ambos habilidades matemáticas que permitam a mediação desses indivíduos em problemas reais da sociedade atual. Por isso, embora apresente definições próprias e diferentes, incluímos aqui a modelagem matemática e a matemática aplicada como atividades de Relações Matemáticas.

## Modelos de Observações

Quando dizemos modelos de observações, nos referimos às oportunidades mais simples que nos rodeiam, a atividades de Modelagem Matemática e/ou Matemática Aplicada, as situações geradoras de debates e reflexões sobre o uso da matemática no dia a dia da sociedade moderna e tecnológica.

Ao entrar em um supermercado, por exemplo, nota-se que a estrutura que sustenta o telhado tem forma piramidal, isso pode parecer sem importância aos olhos dos clientes, mas o professor atento pode explorar esse fato em uma aula de geometria lançando a seguinte questão:

"Por que foi escolhida a forma de pirâmides para sustentação do telhado do supermercado?"

Depois de ouvir as opiniões dos alunos podem-se trabalhar diversos conceitos, como a estrutura triangular das pirâmides e suas funções comparando-as com paralelepípedos, estabelecendo as relações de estabilidade existente nas formas triangulares o que não ocorre nos retângulos, isso é facilmente demonstrado com palitos de sorvetes presos por percevejos de metal. (figura1)

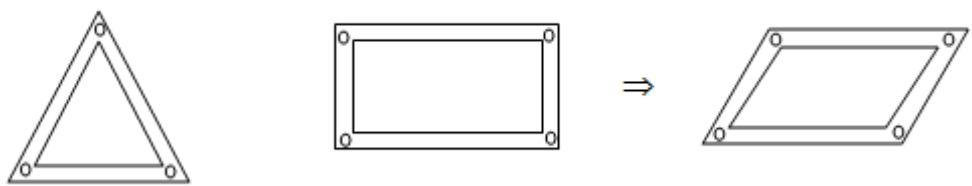

(Figura 1)

Esse simples exemplo ainda pode ser utilizado para explicar, por exemplo, postulados de geometria, explicando que a forma do triângulo é estável por satisfazer o postulado de determinação de plano:

 *"Três pontos não colineares determinam um único plano."*

A observação atenta dos ambientes que nos cercam revela grandes possibilidades de relações matemáticas uteis em sala de aula, cabe ao professor treinar sua capacidade de observação e relação para que perceba inúmeras possibilidades diferentes. Além disso, é imprescindível que o professor de matemática leia, assista filmes, participe de congressos, pois tais hábitos são fontes inesgotáveis de conhecimento das mais diversas áreas. Os filmes e seriados atuais além de divertir estão recheados de informações matemáticas e cálculos que os alunos vêem na escola. O que difere a escola das telas de TV é a roupagem que são apresentados os conceitos, por exemplo, o seriado *SCI: Investigação Criminal,* da TV americana e exibida no Brasil pela Rede Record é assunto diário nos intervalos e nas trocas de professores, até mesmo durante as aulas os comentários sobre situações matemáticas ou científicas são trazidas pelos alunos. O programa está repleto de raciocínio lógico dedutivo, que os alunos se encantam além de trazerem diversos cálculos em seus episódios, outro exemplo de filme é *"O Núcleo: Missão ao Centro da Terra"* da *Paramount Pictures,* uma aventura interessante e cheia de ciência, geografia e matemática. Em determinado momento da trama, por exemplo, as vidas das personagens dependem do cálculo de uma progressão geométrica, eis uma bela oportunidade de interagir com a sala de aula despertando o interesse perdido há algum tempo. O professor atento não precisa de muito esforço para perceber imagens, textos, locais e situações que podem facilmente serem levadas para dentro da sala de aula como tema gerador de suas aulas.

**Considerações Finais**

Boa parte do desinteresse certamente é reflexo de aulas desconexas com a realidade, de conteúdos sem ligações práticas ao dia a dia. As demonstrações mais simples, não carecem de estudos avançados nem de cursos específicos, embora eles sejam extremamente importantes e são de muita valia no aperfeiçoamento do professor de matemática. Oportunidades de aplicações matemáticas surgem naturalmente de observações de situações e ambientes que envolvem o cotidiano de todas as pessoas. O professor atento reúne em si condições suficientes e básicas de estabelecer relações e desenvolver aplicações; é função do professor apresentar a matemática como ferramenta útil à sociedade.

Talvez num primeiro momento, essas aplicações pareçam difíceis de serem percebidas, mas a prática contínua levará ao aperfeiçoamento do professor e das atividades propostas, levando à diminuição do desinteresse, uma vez que estabelecida as relações currículo/realidade, através dos assuntos presentes nas rodas de alunos, estes tendem a compreender a importância da matemática aceitando-a mais facilmente, mesmo que continuem não gostando da disciplina.

# REFERÊNCIAS BIBLIOGRAFICAS